\documentstyle[12pt]{article}
\title{Topological properties defined in terms of generalized
open sets\thanks{1991 Math.\ Subject Classification --- Primary:
54-06, 54D30, 54G12, Secondary: 54A05, 54H05, 54G99.
\protect\newline Key words and phrases --- sg-closed, sg-open,
sg-compact, N-scattered, topological ideal. \protect\newline
Research supported partially by the Ella and Georg Ehrnrooth
Foundation at Merita Bank, Finland and by the Japan - Scandinavia
Sasakawa Foundation.}}
\author{Julian Dontchev\\University of Helsinki\\Department of
Mathematics\\PL 4, Yliopistonkatu 15\\00014 Helsinki\\Finland}
\date{}
\begin{document}
\baselineskip=20pt plus 1pt minus 1pt
\maketitle
\begin{abstract} 
This paper covers some recent progress in the study of sg-open
sets, sg-compact spaces, N-scattered spaces and some related
concepts. A subset $A$ of a topological space $(X,\tau)$ is
called sg-closed if the semi-closure of $A$ is included in every
semi-open superset of $A$. Complements of sg-closed sets are
called sg-open. A topological space $(X,\tau)$ is called
sg-compact if every cover of $X$ by sg-open sets has a finite
subcover. N-scattered space is a topological spaces in which
every nowhere dense subset is scattered.
\end{abstract}

\section{Prelude}

Major part of the talk I presented in August 1997 at the
Topological Conference in Yatsushiro College of Technology is
based on the following three papers:

$\bullet$ {J. Dontchev and H. Maki}, {On sg-closed sets and
semi-$\lambda$-closed sets}, {\em Questions Answers Gen.
Topology}, (Osaka, Japan), {\bf 15} (2) (1997), to appear.

$\bullet$ {J. Dontchev and M. Ganster}, {More on sg-compact
spaces}, {\em Portugal. Math.}, {\bf 55} (1998), to appear.

$\bullet$ {J. Dontchev and D. Rose}, {On spaces whose nowhere
dense subsets are scattered}, {\em Internat. J. Math. Math.
Sci.}, to appear.

The paper in these Proceedings is a collection of the results
from the above mentioned papers and of few new ideas and
questions.

\section{Sg-open sets and sg-compact spaces}

In 1995, sg-compact spaces were introduced independently by
Caldas \cite{CC1} and by Devi, Balachandran and Maki \cite{DBM1}.
A topological space $(X,\tau)$ is called {\em sg-compact} if
every cover of $X$ by sg-open sets has a finite subcover.

Sg-closed and sg-open sets were introduced for the first time by
Bhattacharyya and Lahiri in 1987 \cite{BL1}. Recall that a subset
$A$ of a topological space $(X,\tau)$ is called {\em sg-open}
\cite{BL1} if every semi-closed subset of $A$ is included in the
semi-interior of $A$. A set $A$ is called {\em semi-open} if $A
\subseteq \overline{{\rm Int} A}$ and {\em semi-closed} if ${\rm
Int} \overline{A} \subseteq A$. The {\em semi-interior} of $A$,
denoted by ${\rm sInt}(A)$, is the union of all semi-open subsets
of $A$, while the {\em semi-closure} of $A$, denoted by ${\rm
sCl}(A)$, is the intersection of all semi-closed supersets of
$A$. It is well known that ${\rm sInt}(A)$ = $A \cap
\overline{{\rm Int}A}$ and ${\rm sCl} (A)$ = $A \cup {\rm
Int}\overline {A}$.

Sg-closed sets have been extensively studied in recent years
mainly by (in alphabetical order) Balachandran, Caldas, Devi,
Dontchev, Ganster, Maki, Noiri and Sundaram (see the references).

In the article \cite{BL1}, where sg-closed sets were introduced
for the first time, Bhattacharyya and Lahiri showed that the
union of two sg-closed sets is not in general sg-closed. On its
behalf, this was rather an unexpected result, since most classes
of generalized closed sets are closed under finite unions.
Recently, it was proved \cite[Dontchev; 1997]{JD1} that the class
of sg-closed sets is properly placed between the classes of
semi-closed and semi-preclosed (= $\beta$-closed) sets. All that
inclines to show that the behavior of sg-closed sets is more like
the behavior of semi-open, preopen and semi-preopen sets than the
one of `generalized closed' sets (g-closed, gsp-closed,
$\theta$-closed etc.). Thus, one is more likely to expect that
arbitrary intersection of sg-closed sets is a sg-closed set.
Indeed, in 1997 Dontchev and Maki \cite{DM1} solved the first
problem of Bhattacharyya and Lahiri in the positive.

{\bf Theorem.} \cite[Dontchev and Maki; 1997]{DM1}. An arbitrary
intersection of sg-closed sets is sg-closed.

Every topological space $(X,\tau)$ has a unique decomposition
into two sets $X_1$ and $X_2$, where $X_1 = \{ x \in X \colon \{
x \}$ is nowhere dense$\}$ and $X_2 = \{ x \in X \colon \{ x \}$
is locally dense$\}$. This decomposition is due to Jankovi\'{c}
and Reilly \cite{JR1}. Recall that a set $A$ is said to be {\em
locally dense} \cite[Corson and Michael; 1964]{CM1} (= {\em
preopen}) if  $A \subseteq {\rm Int} \overline{A}$.

It is a fact that a subset $A$ of $X$ is sg-closed (= its
complement is sg-open) if and only if $X_1 \cap {\rm sCl}(A)
\subseteq A$ \cite[Dontchev and Maki; 1997]{DM1}, or equivalently
if and only if $X_1 \cap {\rm Int} \overline{A} \subseteq A$. By
taking complements one easily observes that $A$ is sg-open if and
only if $A \cap X_1 \subseteq {\rm sInt}(A)$. Hence every subset
of $X_2$ is sg-open.

Next we consider the bitopological case and utterly
$(\tau_{i},\tau_{j})$-Baire spaces:

A subset $A$ of a bitopological space $(X,\tau_{1},\tau_{2})$ is
called {\em $(\tau_{i},\tau_{j})$-sg-closed} if $\tau_{j}$-${\rm
Int}(\tau_{i}$-${\rm Cl}(A)) \subseteq U$ whenever $A \subseteq
U$, $U \in SO(X,\tau_{i})$ and $i,j \in \{ 1,2 \}$. Clearly
every $(\tau_{i},\tau_{j})$-rare (= nowhere dense) set is
$(\tau_{i},\tau_{j})$-sg-closed but not vice versa. A subset $A$
of a bitopological space $(X,\tau_{1},\tau_{2})$ is
called {\em $(\tau_{i},\tau_{j})$-rare} \cite[Fukutake,
1992]{F1} if  $\tau_{j}$-${\rm Int}(\tau_{i}$-${\rm Cl}(A))
= \emptyset$, where $i,j \in \{ 1,2 \}$. $A$ is called {\em 
$(\tau_{i},\tau_{j})$-meager} \cite[Fukutake 1992]{F2} if $A$ is
a countable union of $(\tau_{i},\tau_{j})$-rare sets.

A subset $A$ of a bitopological space $(X,\tau_{1},\tau_{2})$ is
called {\em $(\tau_{i},\tau_{j})$-sg-meager} if $A$ is a 
countable union of $(\tau_{i},\tau_{j})$-sg-closed sets. 
Clearly, every $(\tau_{i},\tau_{j})$-meager set is
$(\tau_{i},\tau_{j})$-sg-meager but not vice versa.

{\bf Definition.} \cite[Fukutake; 1992]{F1}. A bitopological
space $(X,\tau_{1},\tau_{2})$ is called {\em
$(\tau_{i},\tau_{j})$-Baire} if $(\tau_{i},\tau_{j})$-${\cal
M} \cap \tau_i = \{ \emptyset \}$, where $i,j \in \{ 1,2 \}$.

{\bf Definition.} A bitopological space $(X,\tau_{1},\tau_{2})$
is called {\em utterly $(\tau_{i},\tau_{j})$-Baire} if
$(\tau_{i},\tau_{j})$-sg-${\cal M} \cap \tau_i = \{ \emptyset
\}$, where $i,j \in \{ 1,2 \}$.

Clearly every utterly $(\tau_{i},\tau_{j})$-Baire space is a
$(\tau_{i},\tau_{j})$-Baire space but not conversely.

{\bf Question 1.} How are utterly $(\tau_{i},\tau_{j})$-Baire
space and $(\tau_{i},\tau_{j})$-semi-Baire spaces related? The
class of $(\tau_{i},\tau_{j})$-semi-Baire spaces was introduced
by Fukutake in 1996 \cite{F2}. Under what conditions is a 
$(\tau_{i},\tau_{j})$-Baire space utterly
$(\tau_{i},\tau_{j})$-Baire?

{\bf Question 2.} Let $(X,\tau_{1},\tau_{2})$ be a bitopological
space such that $\tau_1 \subseteq \tau_2$ and $\tau_2$ is
metrizable and complete. Under what additional conditions is
$(X,\tau_{1},\tau_{2})$ an utterly $(\tau_{i},\tau_{j})$-Baire
space? Note that $(X,\tau_{1},\tau_{2})$ is always a
$(\tau_{i},\tau_{j})$-Baire space \cite[Fukutake; 1992]{F2}.

Observe that sg-open and preopen sets are concepts independent
from each other.

{\bf Theorem.} \cite[Maki, Umehara, Noiri; 1996]{MUN1}. Every
topological space is pre-$T_{\frac{1}{2}}$.

{\bf Theorem.} Every topological space is sg-$T_{\frac{1}{2}}$,
i.e., every singleton is either sg-open or sg-closed.

Improved Jankovi\'{c}-Reilly Decomposition Theorem. Every
topological space $(X,\tau)$ has a unique decomposition into two
sets $X_1$ and $X_2$, where $X_1 = \{ x \in X \colon \{ x \}$ is
nowhere dense$\}$ and $X_2 = \{ x \in X \colon \{ x \}$ is
sg-open and locally dense$\}$.

Let $A$ be a sg-closed subset of a topological space $(X,\tau)$.
If every subset of $A$ is also sg-closed in $(X,\tau)$, then $A$
will be called {\em hereditarily sg-closed} (= hsg-closed)
\cite{DG1}. Hereditarily sg-open sets are defined in a similar
fashion. Observe that every nowhere dense subset is hsg-closed
but not vice versa.

{\bf Theorem.} \cite[Dontchev and Ganster; 1998]{DG1}. For a
subset $A$ of a topological space $(X,\tau)$ the following
conditions are equivalent:

{\rm (1)} $A$ is hsg-closed.

{\rm (2)} $X_1 \cap {\rm Int} \overline{A} = \emptyset$.

A topological space $(X,\tau)$ is called a {\em $C_2$-space}
\cite[Ganster; 1987]{G1} (resp.\ {\em $C_3$-space} \cite{DG1})
if every nowhere dense (resp.\ hsg-closed) set is finite. Clearly
every $C_3$-space is a $C_2$-space. Also, a topological space
$(X,\tau)$ is indiscrete if and only if every subset of $X$ is
hsg-closed (since in that case $X_1 =\emptyset$).

Semi-normal spaces can be characterized via sg-closed sets as
follows:

{\bf Theorem.} \cite[Noiri; 1994]{N1}. A topological space
$(X,\tau)$
is semi-normal if and only for each pair of disjoint semi-closed
sets $A$ and $B$, there exist disjoint sg-open sets $U$ and $V$
such that $A \subseteq U$ and $B \subseteq V$.

{\bf Question 3.} How do hsg-open sets characterize properties
related to semi-normality?

In terms of sg-closed sets, pre sg-continuous functions and pre
sg-closed functions were defined and investigated by Noiri in
1994 \cite{N1}.

Following Hodel \cite{H1}, we say that a {\em cellular family}
in a topological space $(X,\tau)$ is a collection of nonempty,
pairwise disjoint open sets. The following result reveals an
interesting property of $C_2$-spaces.

{\bf Theorem.} \cite[Dontchev and Ganster; 1998]{DG1}. Let
$(X,\tau)$ be a $C_2$-space. Then, every infinite cellular family
has an infinite subfamily whose union is contained in $X_2$.

The {\em $\alpha$-topology} \cite[Nj{\aa}stad; 1965]{Nj1} on a
topological space $(X,\tau)$ is the collection of all sets of the
form $U \setminus N$, where $U \in \tau$ and $N$ is nowhere dense
in $(X,\tau)$. Recall that topological spaces whose
$\alpha$-topologies are hereditarily compact have been shown to
be {\em semi-compact} \cite[Ganster, Jankovi\'{c} and Reilly,
1990]{GJR1}. The original definition of semi-compactness is in
terms of semi-open sets and is due to Dorsett \cite{D1}. By
definition a topological space $(X,\tau)$ is called {\em
semi-compact} if every cover of $X$ by semi-open sets has a
finite subcover.

Remark. (i) The 1-point-compactification of an infinite discrete
space is a $C_2$-space having an infinite cellular family.

(ii) \cite[Ganster; 1987]{G1} A topological space $(X,\tau)$ is
semi-compact if and only if $X$ is a $C_2$-space and every
cellular family is finite.

(iii) \cite[Hanna and Dorsett; 1984]{HD1} Every subspace of a
semi-compact space is semi-compact (as a subspace).

{\bf Theorem.} \cite[Dontchev and Ganster; 1998]{DG1}. {\rm (i)}
Every $C_3$-space $(X,\tau)$ is semi-compact.

{\rm (ii)} Every sg-compact space is semi-compact.

Remark. (i) It is known that sg-open sets are $\beta$-open, i.e.\
they are dense in some regular closed subspace. Note that
$\beta$-compact spaces, i.e.\ the spaces in which every cover by
$\beta$-open sets has a finite subcover are finite \cite[Ganster,
1992]{G2}. However, one can easily find an example of an infinite
sg-compact space -- the real line with the cofinite topology is
such a space.

(ii) In semi-$T_D$-spaces the concepts of sg-compactness and
semi-compactness coincide. Recall that a topological space
$(X,\tau)$ is called a {\em semi-$T_D$-space} \cite[Jankovi\'{c}
and Reilly; 1985]{JR1} if each singleton is either open or
nowhere dense, i.e.\ if every sg-closed set is semi-closed.

{\bf Theorem.} \cite[Dontchev and Ganster; 1998]{DG1}. For a
topological space $(X,\tau)$ the following conditions are
equivalent:

{\rm (1)} $X$ is sg-compact.

{\rm (2)} $X$ is a $C_3$-space.

Remark. (i) If $X_1 = X$, then $(X,\tau)$ is sg-compact if and
only if $(X,\tau)$ is semi-compact. Observe that in this case
sg-closedness and semi-closedness coincide.

(ii) Every infinite set endowed with the cofinite topology is
(hereditarily) sg-compact.

As mentioned before, an arbitrary intersection of sg-closed sets
is also a sg-closed set \cite[Dontchev and Maki; 1997]{DM1}. The
following result provides an answer to the question about the
additivity of sg-closed sets.

{\bf Theorem.} \cite[Dontchev and Ganster; 1998]{DG1}. {\rm (i)}
If $A$ is sg-closed and $B$ is closed, then $A \cup B$ is also
sg-closed.

{\rm (ii)} The intersection of a sg-open and an open set is
always sg-open.

{\rm (iii)} The union of a sg-closed and a semi-closed set need
not be sg-closed, in particular, even finite union of sg-closed
sets need not be sg-closed.

{\bf Problem.} Characterize the spaces, where finite union of
sg-closed sets is sg-closed, i.e.\ the spaces $(X,\tau)$ for
which $SGO(X,\tau)$ is a topology. Note: It is known that the
spaces where $SO(X,\tau)$ is a topology is precisely the class
of extremally disconnected spaces, i.e., the spaces where each
regular open set is regular closed.

A result of Bhattacharyya and Lahiri from 1987 \cite{BL1} states
that if $B \subseteq A \subseteq (X,\tau)$ and $A$ is open and
sg-closed, then $B$ is sg-closed in the subspace $A$ if and only
if $B$ is sg-closed in $X$. Since a subset is regular open if and
only if it is $\alpha$-open and sg-closed \cite[Dontchev and
Przemski; 1996]{DP1}, we obtain the following result:

{\bf Theorem.} \cite[Dontchev and Ganster; 1998]{DG1}. Let $R$
be a regular open subset of a topological space $(X,\tau)$. If
$A \subseteq R$ and $A$ is sg-open in $(R,\tau|R)$, then $A$ is
sg-open in $X$.

Recall that a subset $A$ of a topological space $(X,\tau)$ is
called {\em $\delta$-open} \cite[Veli\v{c}ko; 1968]{V1} if $A$
is a union of regular open sets. The collection of all
$\delta$-open subsets of a topological space $(X,\tau)$ forms the
so called {\em semi-regularization topology}.

{\bf Corollary}. If $A \subseteq B \subseteq (X,\tau)$ such that
$B$ is $\delta$-open in $X$ and $A$ is sg-open in $B$, then $A$
is sg-open in $X$.

{\bf Theorem.} \cite[Dontchev and Ganster; 1998]{DG1}. Every
$\delta$-open subset of a sg-compact space $(X,\tau)$ is
sg-compact, in particular, sg-compactness is hereditary with
respect to regular open sets.

Example. Let $A$ be an infinite set with $p \not\in A$. Let $X
= A \cup \{ p \}$ and $\tau = \{ \emptyset, A, X \}$.

(i) Clearly, $X_1 = \{ p \}$, $X_2 = A$ and for each infinite $B
\subseteq X$, we have $\overline{B} = X$. Hence $X_1 \cap {\rm
Int} \overline{B} \not= \emptyset$, so $B$ is not hsg-closed.
Thus $(X,\tau)$ is a $C_3$-space, so sg-compact. But the open
subspace $A$ is an infinite indiscrete space which is not
sg-compact. This shows that hereditary sg-compactness is a
strictly stronger concept than sg-compactness and '$\delta$-open'
cannot be replaced with 'open'.

(ii) Observe that $X \times X$ contains an infinite nowhere dense
subset, namely  $X \times X \setminus A \times A$. This shows
that even the finite product of two sg-compact spaces need not
be sg-compact, not even a $C_2$-space.

(iii) \cite[Maki, Balachandran and Devi; 1996]{MBD1} If the
nonempty product of two spaces is sg-compact $T_{gs}$-space, then
each factor space is sg-compact.

Recall that a function $f \colon (X,\tau) \rightarrow
(Y,\sigma)$ is called {\em pre-sg-continuous} \cite[Noiri,
1994]{N1} if $f^{-1}(F)$ is sg-closed in $X$ for every
semi-closed subset $F \subseteq Y$.

{\bf Theorem.} \cite[Dontchev and Ganster; 1998]{DG1}. {\rm (i)}
The property 'sg-compact' is topological.

{\rm (ii)} Pre-sg-continuous images of sg-compact spaces are
semi-compact.

\section{N-scattered spaces}

A topological space $(X,\tau)$ is {\em scattered} if every
nonempty subset of $X$ has an isolated point, i.e.\ if $X$ has
no nonempty dense-in-itself subspace. If $\tau^{\alpha} = \tau$,
then $X$ is said to be an {\em $\alpha$-space} \cite[Dontchev and
Rose; 1996]{DR1} or a {\em nodec space}. All submaximal and all
globally disconnected spaces are examples of $\alpha$-spaces.
Recall that a space $X$ is {\em submaximal} if every dense set
is open and {\em globally disconnected} \cite[El'kin; 1969]{El1}
if every set which can be placed between an open set and its
closure is open, i.e.\ if every semi-open set is open.

Recently $\alpha$-scattered spaces (= spaces whose
$\alpha$-topologies are scattered) were considered by Dontchev,
Ganster and Rose \cite{DGR1} and it was proved that a space $X$
is scattered if and only if $X$ is $\alpha$-scattered and
N-scattered.

Recall that a topological ideal $\cal I$, i.e.\ a nonempty
collection of sets of a space $(X,\tau)$ closed under heredity
and finite additivity, is {\em $\tau$-local\/} if $\cal I$
contains all subsets of $X$ which are locally in $\cal I$, where
a subset $A$ is said to be locally in $\cal I$ if it has an open
cover each member of which intersects $A$ in an ideal amount,
i.e.\ each point of $A$ has a neighborhood whose intersection
with $A$ is a member of $\cal I$. This last condition is
equivalent to $A$ being disjoint with $A^{*}({\cal I})$, where
$A^{*}({\cal I}) = \{ x \in X \colon U \cap A \not\in {\cal I}$
for every $U \in \tau_x \}$ with $\tau_x$ being the open
neighborhood system at a point $x \in X$. 

{\bf Definition.} \cite[Dontchev and Rose; 1996]{DR1}. A
topological space $(X,\tau)$ is called {\em N-scattered} if every
nowhere dense subset of $X$ is scattered.

Clearly every scattered and every $\alpha$-space, i.e.\ nodec
space, is N-scattered. In particular, all submaximal spaces are
N-scattered. The density topology on the real line is an example
of an N-scattered space that is not scattered. The space
$(\omega,L)$ below shows that even scattered spaces need not
be $\alpha$-spaces. Another class of spaces that are N-scattered
(but only along with the $T_0$ separation) is Ganster's class of
$C_2$-spaces.

{\bf Theorem.} \cite[Dontchev and Rose; 1996]{DR1}. If $(X,\tau)$
is a $T_1$ dense-in-itself space, then $X$ is N-scattered
$\Leftrightarrow N(\tau) = S(\tau)$, where $N(\tau)$ is the ideal
of nowhere dense subsets of $X$, and $S(\tau)$ is the ideal of
scattered subsets of $X$.

Example. Let $X = \omega$ have the cofinite topology $\tau$. Then
$X$ is a $T_1$ dense-in-itself space with $N(\tau) = I_{\omega}$,
where $I_{\omega}$ is the ideal of all finite sets. Clearly, $X$
is an N-scattered space, since $N(\tau) = I_{\omega} \subseteq
S(\tau)$. Note that $N(\tau) = S(\tau)$. Also, $X$ is far from
being $(\alpha)$-scattered having no isolated points. It may also
be observed that the space of this example is N-scattered being
an $\alpha$-space.

Remark. A space $X$ is called (pointwise) homogenous if for any
pair of points $x,y \in X$, there is a homeomorphism $h \colon
X \rightarrow X$ with $h(x) = y$. Topological groups are such
spaces. Further, such a space is either crowded or discrete. For
if one isolated point exists, then all points are isolated.
However, the above given space $X$ is a crowded homogenous
N-scattered space.

Noticing that scatteredness and $\alpha$-scatteredness are
finitely productive might suggest that N-scatteredness is
finitely productive. But this is not the case.

Example. The usual space of Reals, $(R,\mu)$ is rim-scattered but
not N-scattered. Certainly, the usual base of bounded open
intervals has the property that nonempty boundaries of its
members are scattered. However, the nowhere dense Cantor set is
dense-in-itself. Another example of a rim-scattered space which
is not N-scattered is constructed by Dontchev, Ganster and Rose
\cite{DGR1}.

Remark. It appears that rim-scatteredness is much weaker that
N-scatteredness.

{\bf Theorem.} \cite[Dontchev and Rose; 1996]{DR1}.
N-scatteredness is hereditary.

{\bf Theorem.} \cite[Dontchev and Rose; 1996]{DR1}. The following
are
equivalent:

{\rm (a)} The space $(X,\tau)$ is N-scattered.

{\rm (b)} Every nonempty nowhere dense subspace contains an
isolated point.

{\rm (c)} Every nowhere dense subset is scattered, i.e., $N(\tau)
\subseteq S(\tau)$.

{\rm (d)} Every closed nowhere dense subset is scattered.

{\rm (e)} Every nonempty open subset has a scattered boundary,
i.e., ${\rm Bd} (U) \in S(\tau)$ for each $U \in \tau$.

{\rm (f)} The $\tau^{\alpha}$-boundary of every $\alpha$-open set
is $\tau$-scattered.

{\rm (g)} The boundary of every nonempty semi-open set is
scattered.

{\rm (h)} There is a base for the topology consisting of
N-scattered open subspaces.

{\rm (i)} The space has an open cover of N-scattered subspaces.

{\rm (j)} Every nonempty open subspace is N-scattered.

{\rm (k)} Every nowhere dense subset is $\alpha$-scattered.

{\bf Corollary}. Any union of open N-scattered subspaces of a
space $X$
is an N-scattered subspace of $X$.

Remark. The union of all open N-scattered subsets of a space
$(X,\tau)$ is the largest open N-scattered subset $NS(\tau)$. Its
complement is closed and if nonempty contains a nonempty crowded
nowhere dense set. Moreover, $X$ is N-scattered if and only if
$NS(\tau) = X$. Since partition spaces are precisely those having
no nonempty nowhere dense sets, such spaces are N-scattered. On
the other hand we have the following chain of implications. The
space $X$ is discrete $\Rightarrow$ $X$ is a partition space
$\Rightarrow$ $X$ is zero dimensional $\Rightarrow$ $X$ is
rim-scattered. Also, $X$ is globally disconnected $\Rightarrow$
$X$ is N-scattered. However, this also follows quickly from the
easy to show characterization $X$ is globally disconnected
$\Leftrightarrow$ $X$ is an extremally disconnected
$\alpha$-space, and the fact that every $\alpha$-space is
N-scattered. Actually, something much stronger can be noted.
Every $\alpha$-space is N-closed-and-discrete, i.e.\ $N(\tau)
\subseteq CD(\tau)$. Of course, $CD(\tau) \subseteq D(\tau)
\subseteq S(\tau)$, where $CD(\tau)$ is the ideal of closed and
discrete subsets of $(X,\tau)$, and $D(\tau)$ is the family of
all discrete sets. We will show later that for a non-N-scattered
space $(X,\tau)$ in which $NS(\tau)$ contains a non-discrete
nowhere dense set, $\tau^{\alpha}$ is not the smallest expansion
of $\tau$ for which $X$ is N-scattered, i.e., there exists a
topology $\sigma$ strictly intermediate to $\tau$ and
$\tau^{\alpha}$ such that $NS(\sigma) = X$.

Local N-scatteredness is the same as N-scatteredness.

{\bf Theorem.} \cite[Dontchev and Rose; 1996]{DR1}. If every
point of a space $(X,\tau)$ has an N-scattered neighborhood, then
$X$ itself is N-scattered.

In the absence of separation, a finite union of scattered sets
may fail to be scattered. For example, the singleton subsets of
a two-point indiscrete spaces are scattered. But given two
disjoint scattered subsets, if one has an open neighborhood
disjoint from the other, then their union is scattered.

{\bf Theorem.} \cite[Dontchev and Rose; 1996]{DR1}. In every
$T_0$-space $(X,\tau)$, finite sets are scattered, i.e.,
$I_{\omega} \subseteq S(\tau)$.

{\bf Theorem.} \cite[Dontchev and Rose; 1996]{DR1}. Let
$(X,\tau)$ be a non-N-scattered space, so that $NP(\tau) = X
\setminus NS(\tau) \not= \emptyset$. Suppose also that $NS(\tau)$
contains a nonempty non-discrete nowhere dense subset. Then there
is a topology $\sigma$ with $\tau \subset \sigma \subset
\tau^{\alpha}$ such that $(X,\sigma)$ is N-scattered.

In search for a smallest expansion of $\sigma$ and $\tau$ for
which $(X,\sigma)$ is N-scattered, we have the following:

{\bf Theorem.} \cite[Dontchev and Rose; 1996]{DR1}. Let
$(X,\tau)$ be a space and let $I = \{ A \subseteq E \colon E$ is
a perfect (closed and crowded) nowhere dense subset of
$(X,\tau)\}$. Then $(X,\gamma)$ is N-scattered, where $\gamma =
\tau[I]$, the smallest expansion of $\tau$ for which members of
$I$ are closed.

{\bf Theorem.} \cite[Dontchev and Rose; 1996]{DR1}. Every closed
lower density topological space $(X,F,I,\phi)$ for which $I$ is
a $\sigma$-ideal containing finite subsets of $X$ is N-scattered.
Recall that a lower density space $(X,F,I,\phi)$ is closed if
$\tau_{\phi} \subseteq F$.

{\bf Corollary}. If $(X,F,I,\phi)$ is a closed lower density
space and $I$ is a $\sigma$-ideal with $I_{\omega} \subseteq I$,
then $I = S(\tau_{\phi})$.

{\bf Corollary}. The space of real numbers $R$ with the density
topology $\tau_d$ is N-scattered, and moreover, the scattered
subsets are precisely the Lebesgue null sets.

The following theorem holds and thus we have another (perhaps
new) characterization of $T_0$ separation. A similar
characterization holds for $T_1$ separation.

{\bf Theorem.} \cite[Dontchev and Rose; 1996]{DR1}. A space
$(X,\tau)$ has $T_0$ separation if and only if $I_{\omega}
\subseteq S(\tau)$.

Here is the relation between $C_2$-spaces and $N$-scattered
spaces:

{\bf Corollary}. Every $C_{2}T_{0}$-space is N-scattered.

{\bf Theorem.} \cite[Dontchev and Rose; 1996]{DR1}. A space
$(X,\tau)$ has $T_1$ separation if and only if $I_{\omega}
\subseteq D(\tau)$.

{\bf Theorem.} \cite[Dontchev and Rose; 1996]{DR1}. If $(X,\tau)$
is a $T_0$-space and if $S$ is any scattered subset of $X$ and
if $F$ is any finite subset of $X$, then $S \cup F$ is scattered.

{\bf Corollary}. Every $T_0$-space which is the union of two
scattered subspaces is scattered.

{\bf Corollary}. The family of scattered subsets in a $T_0$-space
is an ideal.

Example. Let $(X,<)$ be any totally ordered set. Then both the
left ray and right ray topologies $L$ and $R$ respectively, are
$T_0$ topologies. They are not $T_1$ if $|X| > 1$. In case $X =
\omega$ with the usual ordinal ordering $<$, $L$ and $R$ are in
fact $T_D$ topologies, i.e.\ singletons are locally closed. The
space $(\omega,L)$, where proper open subsets are finite, is
scattered. For if $\emptyset \not= A \subseteq \omega$ let $n$
be the least element of $A$. Then the open ray $[0,n+1) = [0,n]$
intersects $A$ only at $n$. Thus, $n$ is an isolated point of
$A$. Evidently, $S(L) = P(\omega)$, the maximum ideal. However,
$S(R) = I_{\omega}$, the ideal of finite subsets. For if $A$ is
any infinite subset of $\omega$, $A$ is crowded. For if $m \in
A$ and if $U$ is any right directed ray containing $m$, $(U
\setminus \{ m \}) \cap A \not= \emptyset$, since $U$ omits only
finitely many points of $\omega$. But every finite subset is
scattered. Of course $S(L)$ and $S(R)$ are ideals in the last two
examples since both $L$ and $R$ are $T_D$ topologies.

Remark. Note that the space $(\omega,R)$ is a crowded
$T_D$-space, which is the union of an increasing (countable)
chain of scattered subsets. In particular, $\omega = \cup_{k \in
\omega} \{ n < k \colon k \in \omega \}$ and for each $k$, $\{
n < k \colon k \in \omega \} \in I_{\omega} \subseteq S(R)$. This
seems to indicate that it is not likely that an induction
argument on the cardinality of a scattered set $F$ can be used
similar to the above argument to show that $S \cup F$ is
scattered if $S$ is scattered.

{\bf Question 4.} Characterize the spaces where every hsg-closed
set is scattered. How are they related to other classes of
generalized scattered spaces? Note that:

\begin{center}
Scattered $\Rightarrow$ hsg-scattered $\Rightarrow$ $N$-scattered
\end{center}

Note that the real line with the cofinite topology is an example
of a hsg-scattered space, which is not scattered, while the real
line with the indiscrete topology provides an example of an
$N$-scattered space that is not hsg-scattered.

More generally, if $\cal I$ is a topological (sub)ideal on a
space $(X,\tau)$, investigate the class of $\cal I$-scattered
spaces, i.e.\ the spaces satisfying the condition: ``Every $I \in
{\cal I}$ is a scattered subspace of $(X,\tau)$".

Note that:

$\cal F$-scattered $\Leftrightarrow$ $T_0$-space

$\cal C$-scattered $\Leftrightarrow$ ?

$\cal N$-scattered $\Leftrightarrow$ $N$-scattered space

$\cal M$-scattered $\Leftrightarrow$ ?

${\cal P}(X)$-scattered $\Leftrightarrow$ Scattered space

Of course, every space is ${\cal C}{\cal D}$-scattered.

\
E-mail: {\tt dontchev@cc.helsinki.fi},
{\tt dontchev@e-math.ams.org}
\
\end{document}